# A Stochastic Optimal Control Formulation for Mine Counter Measure Simulations with Multiple Autonomous Survey Vehicles

Philippe Blondeel, Royal Military Academy, Avenue de la Renaissance 30, 1000 Brussels, Belgium, philippe.blondeel@mil.be, ORCID: 0000-0002-2322-095X

Filip Van Utterbeeck, Royal Military Academy, Avenue de la Renaissance 30, 1000 Brussels, Belgium, Filip.VanUtterbeeck@mil.be, ORCID: 0009-0000-5043-6257

Ben Lauwens, Royal Military Academy, Avenue de la Renaissance 30, 1000 Brussels, Belgium, Ben.Lauwens@mil.be, ORCID: 0000-0003-0761-6265

## ABSTRACT

*Modelling and simulating mine counter measure search missions performed by autonomous vehicles equipped with a sensor capable of detecting mines at sea is a challenging endeavour. To address this, we formulated and implemented the problem as a stochastic optimal control model. Our implementation computes an optimal path within a user chosen quadrilateral domain such that the mission duration is minimized for a given residual risk of undetected sea mines. First, we compare the stochastic optimal control implementation against the traditionally used boustrophedon implementation. We show that the mission duration in case of the stochastic optimal control implementation is shorter. Then, by building on our previous work, we introduce a novel mathematical approach that enables multiple autonomous survey vehicles to investigate the domain concurrently. We present results for up to six vehicles, including computed trajectories and an analysis of how mission duration varies with the number of vehicles. Our findings show that mission time decreases non-linearly, , i.e., we observe diminishing returns as more vehicles are added.*

## 1.0   INTRODUCTION

Modelling and simulating mine countermeasures (MCM) search missions performed by autonomous vehicles is a challenging endeavour. The goal of these simulations typically consist of computing an optimal trajectory for the autonomous vehicle in a designated domain where the presence of underwater sea mines is suspected. This trajectory is such that the residual MCM risk in the domain is below a certain threshold, while it is simultaneous ensured that the path time of the autonomous vehicle, i.e., the mission time, is below a certain value. This type of problem is commonly referred to as a Coverage Path Planning (CPP) problem, see [1, 2, 3, 4, 5, 6, 7, 8]. In our previous work [9, 10], we have chosen to model this problem in a stochastic optimal control framework [11], where we minimize the path time for a given residual MCM risk. An important constraint in our implementation consists of the MCM residual risk constraint. This constraint is formulated as an expected value integral over the quadrilateral search domain. In the present work, we first present a comparison between our stochastic optimal control implementation and the traditionally used boustrophedon, i.e., lawnmower, approach. Then, we extend our previous implementation to allow multiple autonomous vehicles to be assigned to the same domain, enabling them to survey the domain in parallel. The structure of this paper is as follows. First, we present our methodology where we give an overview of the governing equations of the sensor attached to the autonomous vehicle. Hereafter, we introduce the mathematical formulation of the stochastic optimal control MCM problem. Then, we extend our mathematical formulation in order to accommodate multiple vehicles working in parallel in a domain. We conclude the methodology part by giving a short description about the numerical implementation of the stochastic optimal control problem. Hereafter, we present the comparison, and our results where we show the trajectories for up to six vehicles inspecting a domain in parallel. Finally, we offer a conclusion.

## 2.0 METHODOLOGY

In this section, we present the methodology. First, we give a short overview of the sensor model. Hereafter, we present the equations governing the stochastic optimal control model for one vehicle. Then, we present our extension as to allow for multiple vehicles surveying the domain concurrently. Last, we give a short description about the numerical implementation.

### 2.1 THE SENSOR MODEL

We briefly introduce the equations pertaining to the modelisation of a Forward Looking Sensor (FLS). A more detailed description can be found in [9, 12].

The sensor model is given by

$$\gamma(x(t), \omega) := \lambda\, p(x(t), \omega)\, F_\alpha(x(t), \omega)\, F_\varepsilon(x(t), \omega),$$

where $\lambda$ stands for the Poisson scan rate in $s^{-1}$.

The sensor model consists of three parts. The first part consists of $p(x(t), \omega)$ and is given by

$$p(x(t), \omega) := \Phi_n\left(\frac{\text{FOM} - 20\log_{10}(\|\omega - x(t)\| + a\|\omega - x(t)\|)}{\sigma}\right)$$

where $\Phi_n(\cdot)$ stands for the cumulative density function (CDF) of the normal distribution, FOM is a parameter related to the sonar characteristics, $a$ is the attenuation coefficient in dB/km, and $x(t)$ is given in as $x(t) := f(x(t), y(t), \psi(t), r(t))$, and $\|\omega - x(t)\| := \sqrt{(\omega_x - x(t))^2} + \sqrt{(\omega_y - y(t))^2}$ where $\omega$ is defined as $\omega := (\omega_x, \omega_y)$, and stands for the position of the possible target, i.e., a sea mine. In $x(t) := f(x(t), y(t), \psi(t), r(t))$, $x(t)$ is the position along the x-axis, $y(t)$ is the position along the y-axis., $\psi(t)$ represents the angle the autonomous vehicle has with respect to the horizontal axis in degrees, and $r(t)$ stands for the turn rate in degrees per seconds. This second part, $F_\alpha(x(t), \omega)$, models the detection of the sensor in front of the autonomous vehicle according to its Field Of View (FOV). The horizontal FOV is schematically shown in Figure 1.

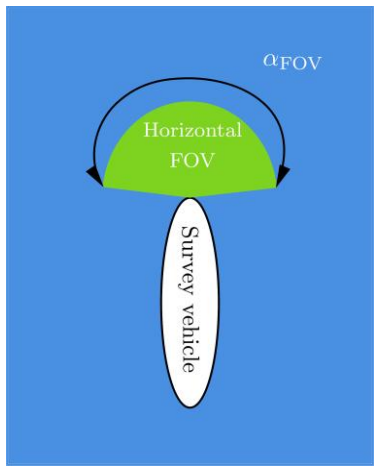

**Figure 1: Top-down view of the autonomous survey vehicle where the Horizontal FOV is shown in green and the water in blue.**

$F_\alpha(x(t), \omega)$ is given as

$$F_\alpha(x(t), \omega) := \frac{1}{1+e^{p_\alpha\left(\frac{-\alpha_{FOV}}{2}-\alpha^b(x(t),\omega)\right)}} + \frac{1}{e^{p_\alpha\left(\alpha^b(x(t),\omega)-\frac{\alpha_{FOV}}{2}\right)}} - 1,$$

where $p_\alpha$ is a parameter used to adjust the slope of the sigmoidal curves, and

$$\alpha^b(x(t), \omega) := \arctan_2\left(dx^b(x(t), \omega), dy^b(x(t), \omega)\right)$$

$$dx^b(x(t), \omega) := (\omega_x - x(t))\cos(\psi(t)) + (\omega_y - y(t))\sin(\psi(t))$$

$$dy^b(x(t), \omega) := -(\omega_x - x(t))\sin(\psi(t)) + (\omega_y - y(t))\cos(\psi(t)),$$

where $\alpha_{FOV}$ is the Field Of View angle of the FLS in degrees. The third part, $F_\varepsilon(x(t), \omega)$ accounts for the height $h$ in meters above the ocean floor of the sensor attached to the autonomous vehicle. A schematic representation is shown in Figure 2.

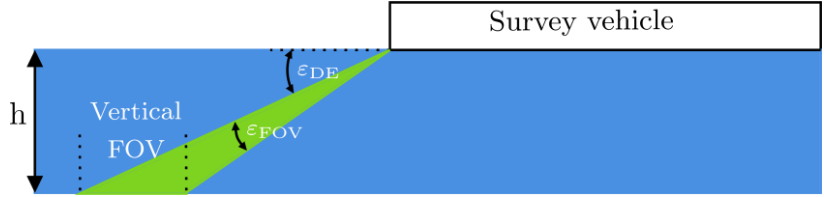

**Figure 2: Schematic side view of the autonomous survey vehicle where the vertical FOV is shown in green and the water in blue.**

$F_\varepsilon(x(t), \omega)$ is given as

$$F_\varepsilon(x(t), \omega) := \frac{1}{1+e^{p_\varepsilon\left(\varepsilon_{DE}-\frac{\varepsilon_{FOV}}{2}-\varepsilon^b(x(t),\omega)\right)}} + \frac{1}{1+e^{p_\varepsilon\left(\varepsilon^b(x(t),\omega)-\varepsilon_{DE}-\frac{\varepsilon_{FOV}}{2}\right)}} - 1,$$

where $p_\varepsilon$ is a parameter used to adjust the slope of the sigmoidal curves, and

$$\varepsilon^b(x(t), \omega) := \arctan\left(\frac{-h}{\|\omega - x(t)\|}\right),$$

where $\varepsilon_{FOV}$ and $\varepsilon_{DE}$ respectively stand for the vertical FOV, and the downward elevation angle such that the sensor can ensonify the sea floor. The final combination of all the different parts of the sensor equation is shown in Figure 3 when the autonomous survey vehicle is present at point (15,15). For a more thorough description we refer to [9, 12].

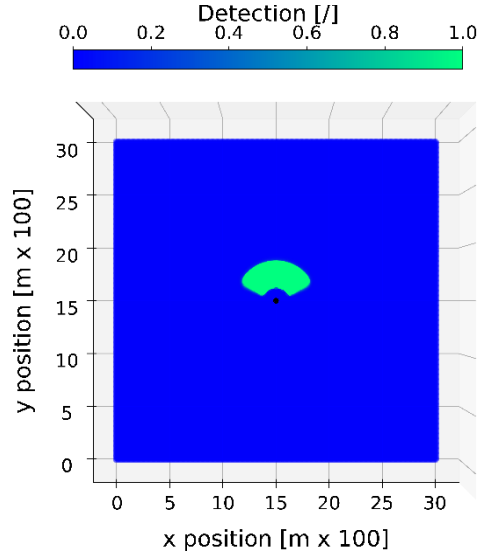

**Figure 3: Illustration of the forward looking detection pattern $\gamma(x(t), \omega)$ where the green area depicts the part that is seen by the onboard sensor. The autonomous survey vehicle is static at location (15,15) (represented by the black dot).**

## 2.2 THE STOCHASTIC OPTIMAL CONTROL FORMULATION

Here and in our previous work [9, 10], the optimization problem is formulated such that the path time of the vehicle (the mission time) $T_F$ needed to survey a designated domain $\Omega$ is minimized for a given user defined residual MCM risk,

$$\min T_F$$

subjected to

$$E[q(T_F)] := \int_\Omega e^{-\int_0^{T_F} \gamma(x(\tau),\omega)d\tau} f_x(\omega)d\omega \leq \text{residual MCM risk}$$

where $\gamma$ represents the sensor function for the FLS and $f_x(\cdot)$ is the probability density function (PDF) of the distribution against which we integrate. In this case we consider $f_x(\cdot)$ to be the PDF of the uniform distribution on a bounded interval $\Omega$,

$$f_x(\omega) = \begin{cases} \frac{1}{|\Omega|}, & \text{if } \omega \in \Omega \\ 0, & \text{otherwise} \end{cases},$$

where $|\Omega|$ stands for the area of the considered domain. This is because we consider the locations of the possible targets, i.e., the sea mines, to be uniformly distributed in a bounded domain. The position of the autonomous vehicle, $x(t)$ at time $t$ is governed by the following differential equations,

$$\frac{dx(t)}{dt} = V\cos(\psi(t))$$
$$\frac{dy(t)}{dt} = V\sin(\psi(t))$$
$$\frac{d\psi(t)}{dt} = r(t)$$
$$\frac{dr(t)}{dt} = \frac{1}{T}\big(Kd(t) - r(t)\big)$$

where $V$ stands for the speed in $m/s$, $K$ is the Nomoto gain constant with units $s^{-1}$, $T$ is the Nomoto time constant with units $s$, and $d(t)$ is the rudder deflection angle given in degrees. The optimization software [11], will control the values for the rudder deflection angle in order to compute a trajectory satisfying the target function and the constraints. The parameters which will be used for all simulations are listed in Table 1.

Table 1: Parameters for the simulations.

| Param. | $\alpha_{FOV}$ | h | $\sigma$ | $\lambda$ | a | $\varepsilon_{FOV}$ | $\varepsilon_{DE}$ | FOM | $p_\alpha$ | $p_\varepsilon$ | V | T | K |
|---|---|---|---|---|---|---|---|---|---|---|---|---|---|
| Value | 120.0 | 20.0 | 9.0 | 20.0 | 5.2 | 5.0 | -6.0 | 72.0 | 25.0 | 400.0 | 2.5 | 0.5 | 5.0 |
| Unit | deg. | $m$ | dB | $s^{-1}$ | db/km | deg. | deg. | / | / | / | $m/s$ | s | $s^{-1}$ |

## 2.3 THE STOCHASTIC OPTIMAL CONTROL FORMULATION FOR MULTIPLE VEHICLES

In order to allow for the computation of trajectories for multiple vehicles operating in parallel in the same domain, the residual risk integral from the previous section is modified as follows for an arbitrary number of vehicles $k$,

$$E[q(T_F)] := \sum_{n=1}^{k} \int_{\Omega} e^{-\int_0^{T_F} \gamma(x_n(\tau),\omega)d\tau} f_x(\omega)d\omega \leq \text{residual MCM risk}$$

where $x_n(t)$ denotes the position of the $n^{th}$ vehicle at the different times $t$. Modifying the residual risk integral also means that the governing differential equations need to be updated. Again for $k$ vehicles, the set of differential equations is modified as follows,

$$n = 1,2,\ldots,k$$
$$\frac{dx_n(t)}{dt} = V_n\cos(\psi_n(t))$$
$$\frac{dy_n(t)}{dt} = V_n\sin(\psi_n(t))$$
$$\frac{d\psi_n(t)}{dt} = r_n(t)$$
$$\frac{dr_n(t)}{dt} = \frac{1}{T_n}\big(K_n d_n(t) - r_n(t)\big)$$

where the underlying assumption is that the speeds of all vehicles is equal, i.e., $\forall n \in \{1,2,\ldots,k\}, V_n = V$, as well as their Nomoto constants, $\forall n \in \{1,2,\ldots,k\}, T_n = T$, and $\forall n \in \{1,2,\ldots,k\}, K_n = K$.

## 2.4 NUMERICAL IMPLEMENTATION

The numerical implementation has been done in the Julia programming language by means of the InfiniteOpt.jl package [11]. Each integral in the residual MCM risk formulation is approximated by a quasi-Monte Carlo integration scheme see [9, 10] as follows

$$\int_\Omega e^{-\int_0^{T_F} \gamma(x(\tau),\omega)d\tau} f_x(\omega)d\omega = \frac{1}{R}\sum_{r=1}^{R} \frac{1}{N}\sum_{n=1}^{N} e^{-\int_0^{T_F} \gamma\left(x(\tau),F^{-1}(u^{(n,r)})\right)d\tau},$$

where $u(n,r)$ is a two-dimensional vector, $u(n,r) := \left[u_1^{(n,r)}, u_2^{(n,r)}\right]$ consisting of shifted qMC quadrature points defined in the unit cube, $F^{-1}(\cdot)$ is the inverse uniform cumulative distribution function applied pointwise to each component of the vector $u^{(n,r)}$ i.e., $F^{-1}(u^{(n,r)}) := \left[F_x^{-1}\left(u_1^{(n,r)}\right), F_x^{-1}\left(u_2^{(n,r)}\right)\right]$, and $F_x^{-1}(\cdot)$ stands for the univariate uniform distribution.

## 3.0 RESULTS

In this section, we first present a comparison between the boustrophedon approach, i.e., the lawnmower approach, and the stochastic optimal control approach in case of one autonomous survey vehicle. Hereafter, we present the results of our stochastic optimal control implementation for multiple autonomous survey vehicles. We are able to simulate up to six individual vehicles operating concurrently in the same domain. All simulations have been performed on a 16 core Intel i7-12850HX processor with 32GB RAM. The different parameters for our simulations are shown in Table 1. We consider a domain $\Omega := [5 \times 25]^2$

### 3.1 THE BOUSTROPHEDON APPROACH COMPARED WITH THE STOCHASTIC OPTIMAL CONTROL APPROACH

In this section we present a comparison between the boustrophedon approach and the stochastic optimal control approach for the case of a single survey vehicle. For the boustrophedon approach, the starting point was chosen as (4.0,7.8) such that, at start, the lower edge of the sensor swath in the y-direction lies entirely, yet just barely, within the survey domain. In this approach, it is common practice for the survey vehicle to perform turns outside of the survey area. As is the case with the stochastic optimal control approach, we require that the survey vehicle returns to its starting point when the residual risk constraint is satisfied. For a residual risk of at most 5%, the trajectories for the boustrophedon and the stochastic optimal control approach are shown in Figure 4 and Figure 5. The path times, i.e., the time it takes for the trajectory to be completed, are given in Table 2.

*Table 2: Path times of the boustrophedon and stochastic optimal control approach.*

|  | Boustrophedon | Stochastic Optimal Control |
|---|---|---|
| Path time [sec]<br>(Mission time) | 4500.00 | 2678.14 |

As can be observed from the path times presented in Table 2, the time needed for the survey vehicle to complete a boustrophedon pattern is longer than in case of the stochastic optimal control formulation. For this specific scenario, we observe a speedup of a factor of 1.68 in favour of the stochastic optimal control approach.

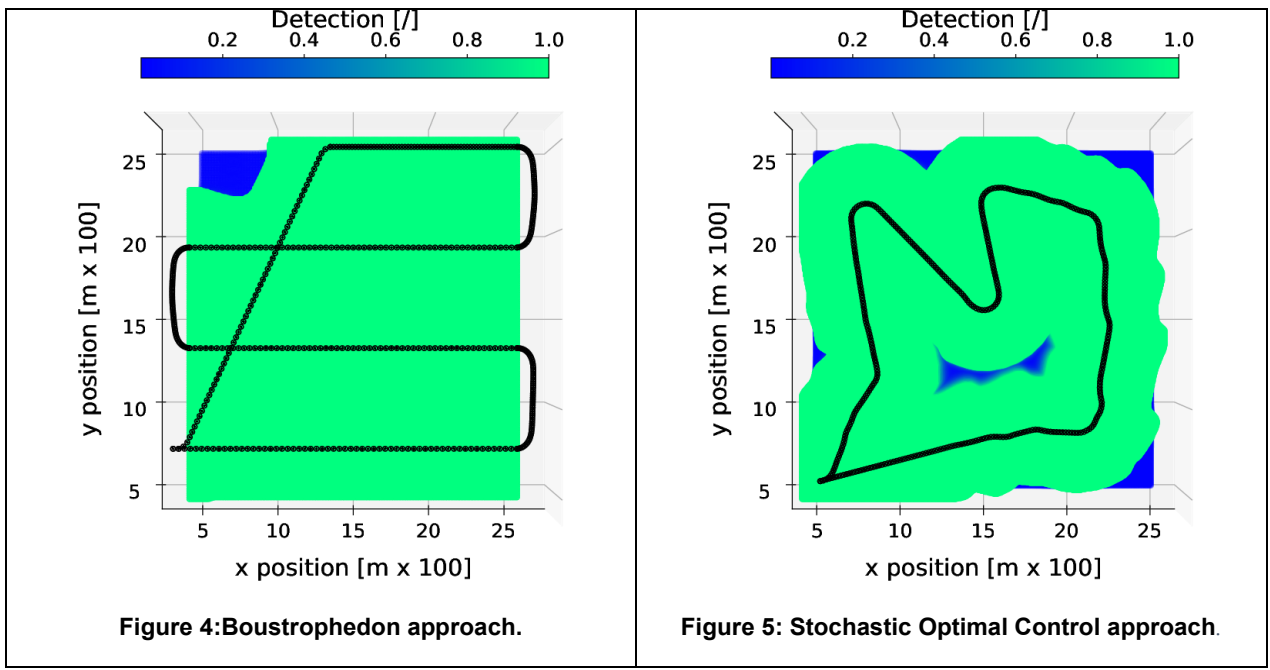

**Figure 4:Boustrophedon approach.**                **Figure 5: Stochastic Optimal Control approach**.

## 3.2 MULTIPLE SURVEY VEHICLES IN THE STOCHASTIC OPTIMAL CONTROL FRAMEWORK

In this section, we present the results of our extended stochastic optimal control implementation, which allows the simulation of up to six vehicles in parallel.

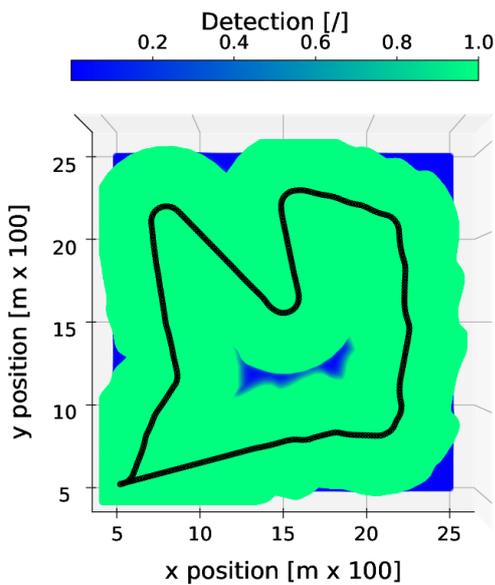    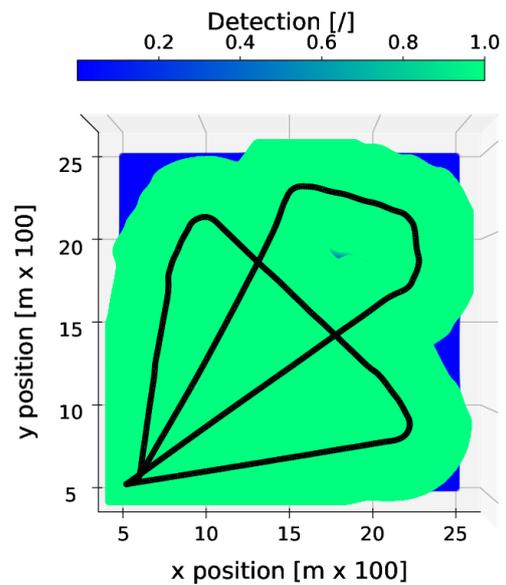

**Figure 6:Trajectory of one vehicle.**                **Figure 7:Trajectories of two vehicles.**

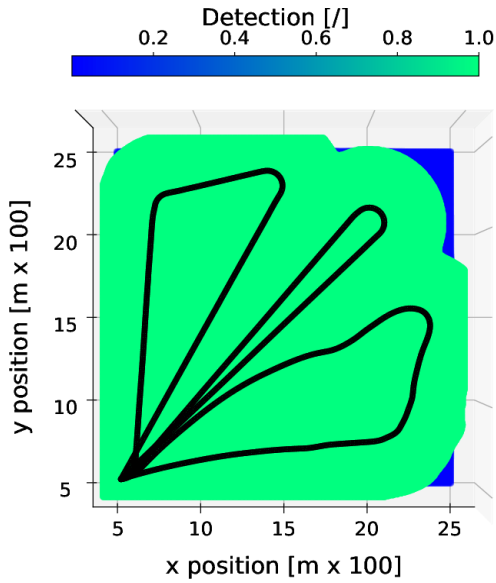

**Figure 8: Trajectories of three vehicles.**

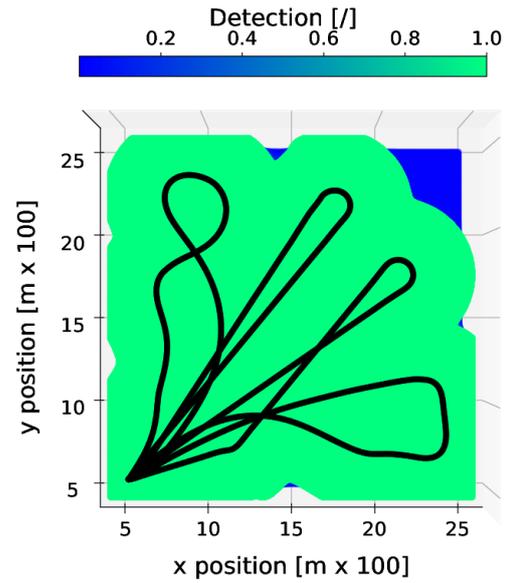

**Figure 9: Trajectories of four vehicles.**

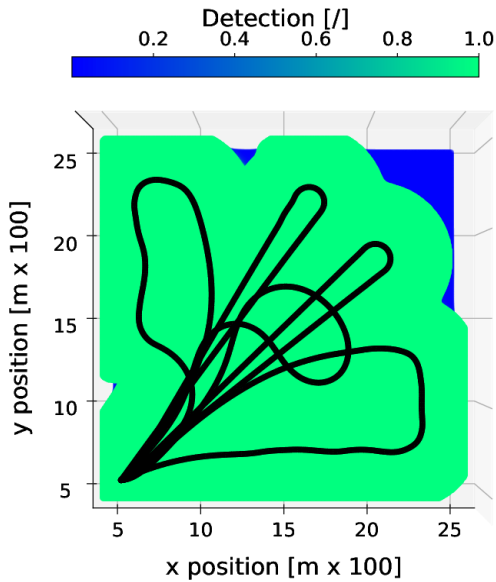

**Figure 10: Trajectories of five vehicles.**

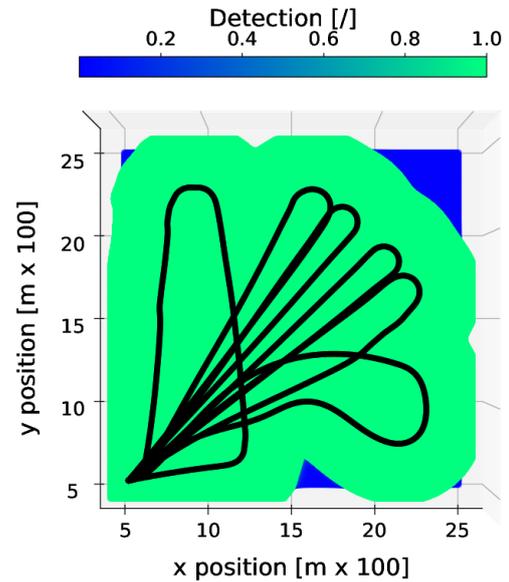

**Figure 11: Trajectories of six vehicles.**

In Figure 6, Figure 7, Figure 8, Figure 9, Figure 10, and Figure 11, we respectively show the trajectories for one, two, three, four, five, and six vehicles surveying the domain in parallel when requesting a residual MCM risk of at most 5%. The black lines depict the trajectories of the individual vehicles, the area in green represents the part of the domain that has been seen by the onboard sonar, and the area in blue stands for the part of the domain that has not been seen.

We also analyse how surveying the domain in parallel with multiple vehicles affects the overall path time (mission time). We assume for this scenario, that all vehicles start from the same location, i.e., point (5.1,5.1)

We present the mission times in function of the number of vehicles in Table 3 and visualize it in Figure 12. As expected, we observe that the mission time decreases when additional vehicles are added to the domain. However, this decrease is non-linear. This means that the time savings decrease with each additional vehicle added to the domain. For example, the decrease in time from one to two vehicles is more significant than from two to three.

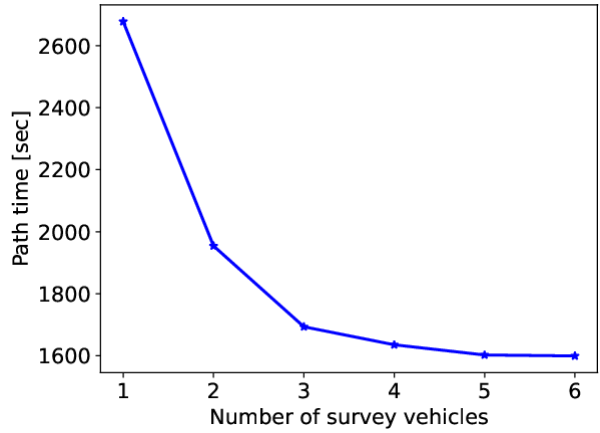

**Figure 12: Path time in function of number of survey vehicles.**

*Table 3: Numerical values of the path time in function of the number of survey vehicles.*

| | Number of survey vehicles | | | | | |
|---|---|---|---|---|---|---|
| | 1 | 2 | 3 | 4 | 5 | 6 |
| Path time [sec] (Mission time) | 2678.14 | 1954.14 | 1693.3 | 1635.17 | 1602.26 | 1599.17 |

## 4.0 CONCLUSION

In this work, we first performed a comparison between the boustrophedon and the stochastic optimal control implementation, where we show that the stochastic optimal control approach yields a speedup of up to a factor 1.68 for the considered scenario. Hereafter, we introduced a novel approach for exploring a domain in parallel using multiple autonomous vehicles, formulated within an optimal control framework. The approach was implemented using the Julia programming language. We have shown the trajectories in a domain for up to six vehicles working in parallel. Additionally, we have investigated in what fashion the mission time decreases when these additional vehicles are added to the domain. We observed a non-linear decrease in time, i.e., adding more vehicles to the domain continues to decrease the time, but this decrease becomes less pronounced with each additional vehicle. In future work, we plan to incorporate the current implementation in a larger framework which will also allow us to compute an optimized schedule of all the different actions of a mine counter measure mission. We plan to base this scheduler on the concept of combinatorial optimization. Furthermore, we plan to incorporate some oceanographic characteristics in our stochastic optimal control simulation such as for example sand ripples which impact the detection capabilities of the sensor.